\documentclass[11pt,twoside]{article}
\usepackage{latexsym}
\usepackage{amssymb,amsbsy,amsmath,amsfonts,amssymb,amscd}
\usepackage{graphicx}
\usepackage{hyperref}
\newcommand{\commentout}[1]{}
\newcommand{\R}{\mathbb{R}}
\newcommand{\N}{\mathbb{N}}

\newcommand {\al} {\alpha}

\newcommand {\sg} {\sigma}

\newcommand {\lb} {\lambda}
\newcommand {\Chi} {{\bf \raise 2pt \hbox{$\chi$}} }
\newcommand {\f}   {\frac}
\newcommand {\p}   {\partial}
\newcommand{\dis}{\displaystyle}
\newcommand {\proof} {\noindent {\bf Proof}. }
\newcommand{\beq}{\begin{equation}}
\newcommand{\eeq}{\end{equation}}
\newcommand{\bea} {\begin{array}{rl}}
\newcommand{\eea} {\end{array}}
\newcommand{\bepa}{\left\{ \begin{array}{ll}}
\newcommand{\eepa} {\end{array}\right.}
\newtheorem{defi}{Definition}
\newtheorem{thm}{Th{e}or{e}m}

\newtheorem{lemma}{Lemma}

\newcommand{\qed}{{ \hfill
                       {\unskip\kern 6pt\penalty 500 \raise -2pt\hbox{\vrule\vbox to 6pt{\hrule width 6pt
                       \vfill\hrule}\vrule} \par}   }}
\title{\Large \bf On a voltage-conductance kinetic system for integrate\&fire neural networks}

\author{
 Beno\^ \i t Perthame\thanks{UPMC Univ Paris 06 and CNRS UMR 7598, Laboratoire Jacques-Louis Lions, F-75005, Paris, France.  Email:~benoit.perthame@upmc.fr}
\thanks{INRIA Paris Rocquencourt. EPC BANG}
\and 
Delphine Salort\thanks{ Institut Jacques Monod UMR 7592  Univ Paris Diderot, Sorbonne Paris Cit{\'e}, F-75205 Paris, France. 
Email: salort.delphine@ijm.univ-paris-diderot.fr}
\and
}

\date{\today}

\begin{document}
\maketitle
\pagestyle{plain}
%\tableofcontents
\pagenumbering{arabic}

\hfill {\em  In memory of late Seiji UKAI, a pioneer in kinetic theory }
\\
\begin{abstract}
The  voltage-conductance kinetic equation  for integrate and fire neurons has been used in neurosciences since a decade and describes the probability density of neurons in a network.  It is used when slow conductance receptors are activated  and noticeable applications to the visual cortex have been worked-out. In  the simplest case, the derivation also uses the assumption of  fully excitatory and moderately  all-to-all coupled networks; this is the situation we consider here.
\\

We study properties of solutions of the kinetic equation for steady states and time evolution and we prove several global a priori bounds both on the probability density and the firing rate of the network. The main difficulties are related to the degeneracy of the diffusion resulting from noise and to the quadratic aspect of the nonlinearity. 
\\

This result constitutes a paradox; the solutions of the  kinetic model, of partially hyperbolic nature, are  globally bounded   but it has been proved that the fully parabolic integrate and fire equation (some kind of diffusion limit of the former) blows-up in finite time. 
\end{abstract}

\bigskip


\bigskip


%%%%%%%%%%%%%%%%%%%%%%%%%%%%%%%%%%%%%%%%%%%%

%%%%%%%%%%%%%%%%%%%%%%%%
%------------------------------------------------------------------
\section{Introduction}
%------------------------------------------------------------------
%%%%%%%%%%%%%%%%%%%%%%%%

Nonlinear Partial Differential Equations arise naturally in the study of neural networks as the closure for a large number of weakly connected neurons in a mean field limit. The complexity of the description of spikes by Hodgkin-Huxley  system has lead many authors to use a simplified yet realistic version called integrate and fire where firing of neurons is assumed when a potential threshold, denoted by $V_F$ below, is achieved. When the network is described solely by the membrane potential $v$ of neurons, the foundations are well established, from physical considerations, comparisons to experimental observations and mathematical theories \cite{BrHa, BrGe,Henry:12, tanre}.  The nonlinearity arises through the total activity of the network (number of spikes per unit of time) at two levels; it generates a current and also an internal noise. For an excitatory network, a recent result is that solutions may blow-up in finite time \cite{CCP}. 
\\

For slow post-synaptic receptors, it is necessary to include the dynamics of conductances  \cite{RKC} and this induces to describe the neurons by the probability density  $p(v,g,t)$ to find neurons  at time $t$ with a membrane potential $v$ and a conductance $g$. This leads to kinetic equations which mathematical structure is reminiscent from the classical Vlasov-Fokker-Planck equation for charged particles \cite{glassey}. A class of such problems arising in neurosciences has been studied  recently by probabilistic methods \cite{PTW1, PTW2}. For networks, the derivation of a mean field equation has been proposed in the last decade for Hodgkin-Huxley or FitzHugh-Nagumo models \cite{bfft}, and for integrate and fire models  \cite{CTSL}. Another direction, still giving rise to a kinetic equation, is to structure the system in voltage and current \cite{BrF}. The models are derived from coupled excitatory neurons where  the membrane potential of each neuron is  governed by a linear integrate and fire equation coupled  with a conductance-based equation which depends on the firing rate of all the neurons.  For integrate and fire models, several studies have been devoted to this kinetic equation as their derivations and applications in particular to the primary visual cortex (V1)  \cite{CCTao, LyT}, their numerical solution \cite{RCT, CCTao}. However, this equation has never been studied from a theoretical mathematical point of view.  In particular the question of blow-up in finite time is open in comparison with the purely parabolic case. This  is our main motivation here; we are going to establish a priori bounds which show that solutions are global. 
\\

These global bounds show a paradox; the kinetic model of hyperbolic nature admits global solutions but the parabolic integrate and fire equation blows-up in finite time. This is similar to the Keller-Segel system which has led to an important literature. The intuition explaining this contradiction is that the times scales in the kinetic and drift-diffusion equations are not the same and the kinetic model is established on a shorter time scale while blow-up occurs at a longer time scale.
\\

The end of this paper is organized as follows. In the next section we present the kinetic model and explain the difficulties to handle it theoretically. 
The third section is devoted to the  stationary state for the linear equation. We prove existence and uniqueness results and some regularity on the solution (see Theorem~\ref{thm:statmain}). In section \ref{sec:ststnl},  we consider the nonlinear kinetic equation; we give some criteria on the strength of synaptic coupling to  classify the regimes where stationary states  exist or not (see Theorem~\ref{thm:fixedpoint}). Section \ref{sec:evol_moments} is concerned with some a priori bounds which ensure that, as long as the solution exists, the propagation of moments holds. Moreover, we obtain some a priori integrability estimates on  the total activity of the network. This allows us to prove in particular that for high  interconnections, the total firing rate cannot be bounded for large time and thus periodic solutions cannot exist. 
The last section is devoted to prove some gain of regularity of the solution and some higher integrability in $L^q$, $q>1$ for the total activity of the network.

%%%%%%%%%%%%%%%%%%%%%%%%
%------------------------------------------------------------------
\section{The voltage-conductance kinetic system}
%------------------------------------------------------------------

The voltage-conductance kinetic system is a nonlinear (2+1) dimensional kinetic Fokker-Plank equation which describes the probability density  $p(v,g,t)$ to find neurons  at time $t$ with a membrane potential $v\in (0, V_F)$ and a conductance $g>0$. It is written as 
\begin{equation}
\label{eq:iftdp}
  \frac{\partial}{\partial t}p(v,g,t)+    \frac{\partial}{\partial v} \left[\big( -g_L v  +g(V_E -v )\big) p(v,g,t)\right] 
+ \frac{\partial}{\partial g} \left[ \f{ g_{\rm in}(t) - g }{\sigma_{E}}  p(v,g,t)\right]  
-  \frac{a(t)}{\sigma_{E}}  \frac{\p^2}{\p g^2} p(v,g,t )  =  0, 
\end{equation}
together with an initial data that satisfies 
$$
p^0(v,g) \geq 0, \qquad \int_0^{V_F} \int_0^\infty p^0(v,g) dv \; dg =1.
$$
The nonlinear aspect comes from the term $g_{\rm in}(t)$. To define it we first  introduce 
\begin{equation}
\label{eq:gfr}
N(g,t):= [- g_L V_F +g(V_E-V_F) ] p(V_F,g,t) \geq 0, \qquad  {\mathcal N}(t) :=  \int_{0}^{+\infty} N(g,t) dg,
\end{equation}
where $N(g,t)$ represents the $g$-dependent firing rate. We
 then define the drift and diffusive coefficients as follows
\begin{equation}\label{defdrift}
g_{\rm in}(t) = f_E \nu(t)+ S_E{\mathcal N}(t),
\end{equation}
\begin{equation}\label{defa}
a(t)= \frac{1}{2 \sigma_{E}} \left( f_E^{2} \nu (t) +   \f{ S_{E}^2}{N_E}  \mathcal{N}(t) \right).
\end{equation}
The parameters that enter this equation have the following interpretations
\\
$\bullet$ $V_E$ is the excitatory reversal potential, 
\\
$\bullet$ Firing occurs when the voltage reaches the threshold  $V_F$,
\\
$\bullet$ Reset is at $V_R$ and we consider that  $0= V_R < V_F< V_E$,
\\
$\bullet$ $g_L > 0$  denotes the leak conductance,
\\
$\bullet$ $g_{\rm in} \geq 0$ is the conductance induced by input currents,
\\
$\bullet$ $\mathcal{N}(t) \geq 0$ is the total firing rate (measures the activity of the network),
\\
$\bullet$ $a(t)= a(\mathcal{N})>0$ represents the intensity of the synaptic noise,
\\
$ \bullet$  $\sigma_E >0$ denotes the time  decay constant of the excitatory conductance,
\\
$\bullet$ $S_{E} \geq 0$ denotes  the synaptic strength of network excitatory coupling,
\\
$\bullet$ $f_{E} > 0$ denotes  the synaptic  strength of the external input $\nu(t)$,
\\
$\bullet$ $N_E $ provides the overall normalization of the coupling strength.
\\

\noindent This kinetic equation  \eqref{eq:iftdp}  is physically  derived in \cite{CTSL}, \cite{RCT} and assumes  that $f_E$ and $\frac{S_E}{f_E}$ are  small enough.

\vspace{0,5cm}

\noindent \underline{Assumptions and notations.}
We assume that the external input rate $\nu(t)$ satisfies for all $t \geq 0$
\begin{equation}\label{asnu}
0< \nu_m  \leq \nu(t) \leq \nu_M < \infty.
\end{equation}
It is sometimes convenient to use notations for the fluxes in  (\ref{eq:iftdp})
\begin{equation}
\label{eq:flux}
J_v (v,g):= ( -g_L v  +g(V_E -v )), \qquad J_g(g) := \sigma_{E}^{-1} ( g_{\rm in}(t) - g ).
\end{equation}
Then, equation (\ref{eq:iftdp}) can be rewritten as
$$  
\frac{\partial}{\partial t}p(v,g,t)+    \frac{\partial}{\partial v} \left[J_vp(v,g,t)\right] 
+ \frac{\partial}{\partial g} \left[  J_g p(v,g,t)\right]  
-  \frac{a(t)}{ \sigma_{E}} \frac{\p^2}{\p g^2} p(v,g,t )  =  0.
$$

\noindent \underline{Boundary conditions.} We need to impose boundary conditions at $v=0,V_F$ and $g=0,+\infty$.

\noindent  To define the boundary conditions on $v=0$ and $v=V_F$, it is assumed that when a neuron reach the threshold voltage $V_F$,  the voltage instantaneously  resets to the value $V_R=0$, without refractory states and that the  conductance stays with the same value upon voltage reset. More precisely, if we set   $g_F= \frac{g_L V_F}{V_E-V_F}$, then, for $g\leq g_F$, the  two  fields $J_v$ at $v=0, \; V_F$ are in-coming and we use the Dirichlet condition 
\begin{equation}
\label{eq:BC1}
p(0,g,t) =p(V_F,g,t)  = 0, \qquad \text{ for  } \;  g  \leq g_F.
\end{equation}
For $g>g_F$, the field $J_v$  is out-going at $v=V_F$ and in-going at $v=0$. The model expresses that neurons undergo a discharge  at $V_F$ and are reset at $V_R=0$ which leads to write 
\begin{equation}
\label{eq:BC2}
g V_E p(0,g,t)= [- g_L V_F +g(V_E-V_F) ] p(V_F,g,t)  \qquad \text{ for } \; g > g_F .
\end{equation}
Notice however that this flux equality holds globally for all $g \in (0, \infty)$. 
\newline
The boundary conditions at $g=0$ and $g=+\infty$ are simply zero flux conditions 
\begin{equation} \label{eq:NF}
(-g+g_{\rm in}) p(v,g,t ) -a \frac{\p}{\p g} p(v,g,t ) =0,    \qquad \text{ for  } \;  g =0,\; g= +\infty.
\end{equation}
Those boundary conditions,  when integrating the  equation \eqref{eq:iftdp}, imply the conservation property
\begin{equation} \label{eq:proba}
\int_0^{V_F} \int_0^\infty p(v,g,t) dv \; dg =1
\end{equation}
which is in accordance with the interpretation that  the solution is a probability density (when the initial data is). 
\\

%------------------------------------------------------------------------
\noindent \underline{Main difficulties.} 
%------------------------------------------------------------------------
% 
The kinetic equation \eqref{eq:iftdp} generalizes the Fokker-Planck-Kolmogorov  equations for network integrate and fire models (see \cite{BrHa,CCP, tanre} and references therein) which describe the dynamic on the neurons only via their potential membrane.  Hence, the theoretical study of equation  \eqref{eq:iftdp} is a priori more complicate than the  Fokker-Planck equations. In particular, all the difficulties encountered for the former arise also in our context in particular the possible blow-up in finite time.

A first difficulty is that the  operator 
$$  
\frac{\partial}{\partial v} \left[\big( -g_L v  +g(V_E -v )\big) p(v,g,t)\right] 
+\frac{1}{\sigma_{E}} \left(\frac{\partial}{\partial g} \left[ ( g_{\rm in}(t) - g ) p(v,g,t)\right]  
- a(t) \frac{\p^2}{\p g^2} p(v,g,t) \right)
$$
has a partial diffusion on the variable $g$ only. A consequence is the difficulty to prove regularity for the solution, even in the linear case ($S_E=0$) or for the stationary states ($g_{\rm in}(t)=cst$).  However, the hypoelliptic character (\cite{Villani_hypo}) of the above operator added to the specificity of dimension 2+1 gives us an opening to obtain some regularity. More precisely, for the  stationary states (when they exist), we  prove some smoothness,  in terms of derivatives, of the solution in  both directions, which implies that the stationary states live in some $L^q$ space, $q>1$. For the evolution equation,  the strategy is completely different, and we  obtain some a priori estimates where  the solution propagates the moments well  in $L^q$, $q\geq 1$ (without proving smoothness in term of derivatives), by using multipliers adapted to the above  operator (see sections \ref{sec:stst} and \ref{sec:reg}).
\newline
The second difficulty comes from the nonlinearity. It is driven by the average  $\mathcal{N}(t)$ of the  boundary flux $N(t,g)$, and it is difficult to prove bounds on these quantities. We obtain a priori estimates where $\mathcal{N}$ is locally in $L^q$, $q>1$ in time (see Theorem~\ref{thml2}), assuming the initial data is sufficiently decreasing at infinity (typically a Gaussian in the variable $g$).  This control of the firing rate points out the difference between  the dynamic of  the kinetic equation and the integrate and fire equation structured only by the potential, where blow-up of solutions arises, even for smooth initial data and weak interconnections. This arises with an $L^1_{\rm loc}$  bound on $\mathcal{N}(t)$ because a Dirac mass forms in finite time, see \cite{CCP}. Moreover, our a priori estimates on $\mathcal{N}(t)$ proved in Theorem~\ref{thml2} are a support to rigorously prove global existence of solutions of equation \eqref{eq:iftdp}, following the Lions-Aubin time compactness argument.  Finally, let us mention that  the very intuitive result of convergence of the solution to the stationary states for weak connections (numerically observed, see \cite{CCTao} and references therein) seems theoretically a difficult question because of the degeneracy of the kinetic equation; in particular, convergence to the stationary state in the case of weak interconnections can not be directly  deduced  by a perturbation argument from our result in the linear case  (see Theorem~\ref{thm:linear}). 
%
%
%%%%%%%%%%%%%%%%%%%%%%%%%%%%%%%%%%%%
%--------------------------------------------------------------------------------------------------------
\section{The steady state for the linear equation ($S_E=0$)}
\label{sec:stst}
%--------------------------------------------------------------------------------------------------------
%%%%%%%%%%%%%%%%%%%%%%%%%%%%%%%%%%%%

Our first goal is to  study  the stationary state associated to the linear equation  \eqref{eq:iftdp}, that is when $g_{\rm in}$ and $a$ are two positive constants. The equation is then given by
\begin{equation} \begin{cases} 
  \frac{\partial}{\partial v} \left[\big( -g_L v  +g(V_E -v )\big) p(v,g)\right] 
+ \frac{\partial}{\partial g} \left[\f{ g_{\rm in} - g } {\sigma_E} \;  p(v,g)\right]  
-    \frac{a}{\sigma_E}  \frac{\p^2}{\p g^2} p(v,g)  =  0,
\\ \\
\displaystyle \int_0^\infty \int_0^{V_F} p(v,g) dv dg=1, \qquad p (v,g) \geq 0, 
\end{cases} 
\label{eq:iftd}
\end{equation}
with the boundary conditions\eqref{eq:BC1}, \eqref{eq:BC2},  \eqref{eq:NF}.
\\

The main difficulty is to prove some a priori regularity estimates on the solution in order to obtain that the solution has a higher integrability than merely $L^1$. Indeed, only partial regularity stems from the diffusion in $g$ but estimates in $v$ come from  the hypoelliptic character of this kinetic equation, sharing thus many similarities with the kinetic Fokker-Planck equations which have been recently  studied by many specialists \cite{bouchut, Villani_hypo, ArseSR}. This structure will allow us to derive some Sobolev regularity and, by variants around Sobolev injections, higher integrability than $L^1$. However we do make use explicitly of commutators in our approach.
These estimates are the first step to prove existence of a nontrivial and nonnegative solution of equation \eqref{eq:iftd} by an approximation argument. 
\\

Our approach  uses the Besov spaces, variants of the classical Sobolev  spaces $W^{s,p}$, defined as follows (see for example \cite{albert}, \cite{bcd}  and references therein)
%
%-------------------------------
\begin{defi}
Let  $\Omega$ be a domain and let $p \geq 1$, $s \geq 0$. We set $\Omega_h:= \{x \in \Omega, x+h \in \Omega \}$ and we define the Besov space
$$
B^{s,p}_{\infty}(\Omega) :=  \Big \{ f \in L^p(\Omega); \;  |f|_{B^{s,p}_\infty(\Omega) }:= \sup_{h \leq 1}\;  h^{-s}\;  \| f(\cdotp+h)-f(\cdotp)\|_{L^{p}(\Omega_h)} < \infty \Big\} . 
$$
Here $s$ can be interpreted as the regularity exponent (number of derivatives) and $p$ as the $L^p$ space to measure the smoothness of the function under consideration.
\end{defi}
%------------------------------- 
We can now state our  main theorem 
%-------------------------------------------------------  
\begin{thm}\label{thm:statmain}
There exists a unique  nonnegative solution $p$ of equation  \eqref{eq:iftd} with the following regularities
\begin{equation}\label{regest1}
  \big| \big( -g_L v  +g(V_E -v )\big) p(v,g)\big|_{B^{1/2, 1 }_{\infty}( (0,V_F)\times \R^+) }<+\infty, 
\end{equation}
\begin{equation}\label{regest2}
   \|p\|_{L^q((0,V_F)\times \R^+)}<+\infty, \qquad  \quad  1 \leq  q  < \frac{8}{7} \cdotp
\end{equation}
\end{thm}
%-------------------------------------------------------  

Let us briefly explain the idea of the proof of Theorem~\ref{thm:statmain}.
To prove the a priori estimates  (\ref{regest1}) and (\ref{regest2}),  we first use the  diffusive part of equation  \eqref{eq:iftd} and basics estimates on the solution. This allows us, on the one hand to gain some regularity in the variable $g$ for the function $p$ and   to prove some regularity  with respect to the variable $v$ of the flux $( -g_L v  +g(V_E -v )) p(v,g)$, but at the cost of low regularity in the variable $g$.
In the second step, we prove  (\ref{regest1}) by an interpolation argument. The third step is devoted to the proof of  (\ref{regest2}). The difficulty is that the weight $J_v= -g_L v  +g(V_E -v )$ vanishes and so we cannot directly apply the Sobolev injections.
 Finally, we prove existence and uniqueness of the linear equation of \eqref{eq:iftd} using the a priori estimates    (\ref{regest1}), (\ref{regest2}).

%%%%%%%%%%%%%%%%%%%%%%%%%%%%%%%
%-------------------------------------------------------------------------------------------
\subsection{Basic estimates on the stationary solution}
%-------------------------------------------------------------------------------------------

Elementary manipulations give several useful basic estimates on solutions of equation \eqref{eq:iftd}. These are, for some constants $Z, \; K_1$ and $K_2$ which we  estimate later on, 
\begin{equation} \label{eq:intv}
\int_0^{V_F} p(v,g) dv = Z(g_{\rm in})^{-1} \; e^{- \frac{1}{a}( g- g_{\rm in} )^2/2} , \qquad \forall g \geq 0,  \hbox{ with } \sqrt{ 2 a\pi }Ê\geq Z(g_{\rm in})  \geq  \sqrt{\f{a\pi}{2}}, 
\end{equation}
\begin{equation}\label{eq:moments}
\int_0^\infty \int_0^{V_F} e^{\frac{  g^2}{4a}}p(v,g) dv dg \leq K_1(g_{\rm in}),
\end{equation}
\begin{equation}\label{est:tfr}
0 \leq \int_{0}^{+\infty} [ -g_L v  +g(V_E -v ) ] \; p(v,g)dg : =  {\mathcal N}  \leq K_2(g_{\rm in}),   \qquad \forall v \in [ 0, V_F] ,
\end{equation}
\begin{equation}\label{estgin}
 g_{\rm in} + \sqrt{ \frac{a }{2\pi}} e^{- \frac{ g_{\rm in}^{2}}{2a} } \leq  \int_0^\infty \int_0^{V_F} g p(v,g) dv dg \leq  g_{\rm in} + \sqrt{ \frac{2 a }{\pi}} e^{- \frac{ g_{\rm in}^{2}}{2a} }.
\end{equation}
\

The first estimate follows after integration of~\eqref{eq:iftd}  in $v$ which gives thanks to the zero flux conditions~\eqref{eq:BC1} and \eqref{eq:BC2} that
$$
 \f{d}{dg} \left( \big( g_{\rm in} - g \big) \int_0^{V_F} p(v,g) dv  -  a\frac{d}{d g} \int_0^{V_F} p(v,g) dv \right) =  0
$$
and so, using the boundary condition  \eqref{eq:NF}, we obtain that
$$ \big( g_{\rm in} - g \big) \int_0^{V_F} p(v,g) dv  - a\frac{d}{d g} \int_0^{V_F} p(v,g) dv =0.$$
Hence, the solution of the above equation is a Maxwellian  with (because $g_{\rm in} \geq 0$) 
$$
\sqrt{ 2 \pi a }Ê\geq Z(g_{\rm in}) := \int_0^{\infty} e^{- \frac{(g - g_{\rm in} )^2}{2 a}}  dg \geq Z(0)= \sqrt{\f{a \pi}{2}}  \cdotp
$$ 
The inequality \eqref{eq:moments} is a direct consequence of  (\ref{eq:intv}).
\\

 The inequality \eqref{estgin} is just a consequence of the integration in $g$, with weight $g$, of equation in  \eqref{eq:intv} whic gives
$$\begin{array}{rl}
 \int_0^\infty \int_0^{V_F} g p(v,g) dv dg & =  g_{\rm in} + Z(g_{\rm in})^{-1}  \int_0^\infty (g-g_{\rm in})e^{- \frac{( g-g_{\rm in} )^2}{2 a}} dg
 \\[3mm]
 & =  g_{\rm in} + a Z(g_{\rm in})^{-1} e^{- \frac{ g_{\rm in}^{2}}{2a} }.
\end{array} $$
It remains to use the bound  $Z(g_{\rm in})$ in \ref{eq:intv}.
\\

For the  estimate \eqref{est:tfr}, we integrate in $g\in(0,\infty)$ equation \eqref{eq:iftd} to obtain 
$$
\f d{dv} \left[ \int_{0}^{+\infty} ( -g_L v  +g(V_E -v ))p(v,g)dg \right]=0,
$$
$$
 \int_{0}^{+\infty} ( -g_L v  +g(V_E -v ))p(v,g)dg=  {\mathcal N} \qquad \forall v \in [0, V_F].
$$
 Integrating this equation between $0$ and $V_F$ and using estimate (\ref{eq:moments}), we obtain \eqref{est:tfr}. 

%%%%%%%%%%%%%%%%%%%%%%%%%%%%%%%
%-------------------------------------------------------------------------------------------
\subsection{Gradient  estimates in $g$}
%-------------------------------------------------------------------
%%%%%%%%%%%%%%%%%%%%%%%

A second family of estimates follows from the diffusion in $g$. We are going to prove the  
%
%-------------------------------------------------------------------
\begin{lemma}\label{lem:estdif} 
The a priori estimates hold for solutions of the  equation \eqref{eq:iftd} 
\beq \label{est:log1}
\int_0^{V_F} \int_0^\infty \big|\nabla_g \sqrt{ p(v,g) } \big|^2 dv\; dg \leq K_3(g_{\rm in}),
\eeq
with
\beq \label{eq:K3}
K_3(g_{\rm in})= Z^{-1}( g_{\rm in}) \frac{1}{2a \sigma_E} \int_{0}^{+\infty}   \big( g_{\rm in} - g \big)^{2}  e^{-\frac{(g-g_{\rm in})^2}{2a}}dg  +g_L+ g_{\rm in} + aZ(g_{\rm in})^{-1} 
 e^{-  \frac{g_{ in}^{2}}{2a} }<+\infty.
\eeq
Consequently,  there is a constant $C$ such that
\beq \label{est:dv}
\int_0^{V_F} \int_0^\infty e^{\frac{g^2}{8a}} \left| \nabla_g p\right| dv dg \leq C , \qquad \int_0^{V_F} \int_0^\infty \left| \f{d}{dv }\left( \int_{g}^{+\infty} p(v,g')dg'\right) \right| dv\; dg \leq C .
\eeq
\label{lm:est1}
\end{lemma}
%-------------------------------------------------------------------
%
{\bf Proof of Lemma~\ref{lem:estdif}.}
We multiply the equation \eqref{eq:iftd} by $\ln p$ and integrate. We find
$$\bea
- \dis \int_0^{V_F} \int_0^\infty  \left[\big( -g_L v  +g(V_E -v )\big) \f{\p p(v,g)}{\p v} + \frac{1}{\sigma_E} \left( \big( g_{\rm in} - g \big) \f{\p p(v,g)}{\p g} -a \f{|\nabla_g p|^2}{p}
 \right)\right] &
\\[10pt]
+ \dis \int_0^\infty \big( -g_L v  +g(V_E -v )\big)  p(v,g) \ln p \Big|_0^{V_F} \; dg &=0 .
\eea$$
We may integrate by parts again and find
$$\bea
\dis  \int_0^{V_F} \int_0^\infty  \left[- \big( g_L +g \big) p(v,g)  + \frac{1}{\sigma_E}a \f{|\nabla_g p|^2}{p}\right] 
+  \frac{1}{\sigma_E}\int_0^\infty \big( -g_L v  +g(V_E -v )\big)  p(v,g) \ln p (v,g) \Big|_0^{V_F} \; dg
\\[5pt]
- \dis  \int_0^\infty  \big( -g_L v  +g(V_E -v )\big)  p(v,g) \Big|_0^{V_F} \; dg - \frac{1}{\sigma_E} \int_{0}^{+\infty}\int_0^{V_F}   \big( g_{\rm in} - g \big)  \partial_gp(v,g) \; dgdv= 0 .
\eea $$
Using \eqref{eq:moments} and  \eqref{est:tfr}, with   $Z$ defined in equation (\ref{eq:intv}), this is reduced to  
\beq
 \frac{a}{\sigma_E}\int_0^{V_F} \int_0^\infty  \f{|\nabla_g p(v,g)  |^2}{p(v,g) } dv\; dg +  \int_0^\infty \big( -g_L v  +g(V_E -v )\big)  p(v,g) \ln p \Big|_0^{V_F} \; dg = K_3(g_{\rm in}).
\label{apslog1}
\eeq
with $K_3(g_{\rm in})$ given by \eqref{eq:K3}.

\vspace{0,5cm}

\noindent On the other hand, we have
\beq \label{est:log2}
0\leq   \int_0^\infty \big( -g_L v  +g(V_E -v )\big)  p(v,g) \ln p \Big|_0^{V_F} \; dg = \ln \f{g V_E }{ gV_E-gV_F- g_L V_F } \; \int_{ g_{F}} ^\infty N(g)  dg \leq  K_3(g_{\rm in}).
\eeq
Indeed, we notice that 
$$
 \int_0^\infty \big( -g_L v  +g(V_E -v )\big)  p(v,g) \ln p \Big|_0^{V_F} \; dg =  \int_{ g_{F}} ^\infty N(g)  \ln \f{p(V_F,g)}{p(0,g)} \; dg
$$
and the zero-flux condition \eqref{eq:BC2} can be written as 
$$
g V_E = [- g_L V_F +g(V_E-V_F) ]\;  \f{p(V_F,g)}{p(0,g)} 
$$
so that 
$$
 \int_0^\infty \big( -g_L v  +g(V_E -v )\big)  p(v,g) \ln p \big|_0^{V_F} \; dg = \int_{ g_{F}} ^\infty N(g) \ln \f{g V_E }{ gV_E-gV_F- g_L V_F } \; dg \geq 0 
$$ 
and (\ref{est:log1}) and (\ref{est:log2}) follow from \eqref{apslog1}.
\\

The first estimate in (\ref{est:dv}) is just a combination of (\ref{est:log1}) writing $ \nabla_g p= 2 \sqrt p  \nabla_g \sqrt p$ and using Cauchy-Schwarz inequality 
$$
\left( \int_{0}^{V_F}\int_{0}^{+\infty} e^{\frac{g^2}{8a}}  | \nabla_g p| dv dg \right)^2 \le 4 \left( \int_{0}^{V_F}\int_{0}^{+\infty}| \nabla_g  \sqrt p|^2 dv dg  \right) \left( \int_{0}^{V_F}\int_{0}^{+\infty}   e^{\frac{g^2}{4a}}  p  dv dg \right).
$$
For the second, let $g \in \R^+$. We integrate equation  (\ref{eq:iftd})  between $g$ and $+\infty$ and, using the zero flux condition at $g=+\infty$, we obtain
$$
 \left| \f d{dv} \left( \int_{g}^{+\infty} ( -g_L v  +g(V_E -v ))p(v,g')dg'\right) \right| \leq \frac{1}{\sigma_E} \left( | g_{\rm in} - g | p(v,g)+ a | \partial_{g}p| \right),
 $$
which completes the proof of estimate (\ref{est:dv}) thanks to the integrability of the right hand side which is already proved, and Lemma~\ref{lem:estdif} follows.
\hfill $\square$

%%%%%%%%%%%%%%%%%%%%%%%%%%%%%%%
%-------------------------------------------------------------------------------------------
\subsection{Besov regularity}
%-------------------------------------------------------------------------------------------
%%%%%%%%%%%%%%%%%%%%%%%%%%%%%%%

We choose $\Omega=(0,V_F)\times (0,+\infty)$ and let $h:=(h_1,h_2)$. To prove the Besov estimate  (\ref{regest1}), it is enough to show that there exists a constant $C$ independent of $h$ such that
$$
\int_{\Omega_h} | \phi(v+h_1,g+h_2) -  \phi(v,g)|dvdg \leq C |h|^{\frac{1}{2}}, \qquad \phi :=  (- g_L v +g(V_E-v))p . 
$$
 We use an interpolation method. The  gain of regularity in the variable $g$ given by estimate (\ref{est:log1})
 compensates the fact that the second estimate in (\ref{est:dv}) gives a gain of regularity in $v$ only with loss of one derivative in the variable $g$. 
\\

We begin with the translate in $g$.  We first notice the following Lemma which is a simple consequence of \eqref{eq:moments} and \eqref{est:dv} and thus we do not prove it, 
%-------------------------------------------------
\begin{lemma}\label{lem:dgflux}
Let $p$ be the solution of equation (\ref{eq:iftd}).  Then
\begin{equation}\label{estdg}
\phi(v,g) =  (-g_L v+g(V_E-v))p \in L^1\big((0,V_F); W^{1,1}(\R^+) \big).
\end{equation}
\end{lemma}
%-------------------------------------------------

From estimates (\ref{estdg}), we deduce that there exists a constant  $C$ such that
\begin{equation}\label{estdgnum}
\int_{0}^{V_F} \int_{0}^{+\infty} \mathbb{I}_{g+h_2 \in (0,+\infty)}| \phi(v,g+h_2)-\phi(v,g) | dvdg \leq C|h_2|
\end{equation}
and thus
$$
\int_{\Omega_h}|\phi(v+h_1,g+h_2)- \phi(v+h_1,g)| dvdg \leq C|h_2| \leq C |h|^{\frac{1}{2}}.
$$
We write
$$
\phi(v+h_1,g+h_2)-\phi(v,g)= \phi(v+h_1,g+h_2)- \phi(v+h_1,g)+\phi(v+h_1,g)-\phi(v,g).
$$
So we are reduced to prove that
\beq
\int_{\Omega_h}|\phi(v+h_1,g)-\phi(v,g)| dvdg  \leq C |h|^{\frac{1}{2}}.
\label{eq:vtrans}
\eeq

To estimate these $v$-translates, we first recall from  (\ref{est:dv}) that there exists a constant  $C$ such that 
\begin{equation}\label{estdvnum}
\int_{0}^{V_F} \int_{0}^{+\infty} \mathbb{I}_{v+h_1 \in (0,V_F)} \left | \int_{g}^{+\infty} \phi(v+h_1,w)-\phi(v,w) dw \right| dvdg \leq C|h_1| .
\end{equation}
 We set $H=|h_1|^{\frac{1}{2}}$ and write 
$$
\phi(v+h_1,g)-\phi(v,g)= \frac{1}{H} \int_{0}^{H} \big[ \phi(v+h_1,g)-\phi(v+h_1,g+s)+ \phi(v+h_1,g+s)-\phi(v,g+s)+\phi(v,g+s)-\phi(v,g) \big]ds
$$
which gives the control
$$
|\phi(v+h_1,g)-\phi(v,g)| \leq A(v,g)+B(v,g)
$$
with
$$
A(v,g):=\frac{1}{H} \int_{0}^{H} |\phi(v + h_1, g)-\phi(v + h_1, g + s)|ds +\frac{1}{H} \int_{0}^{H}	|-\phi(v, g) + \phi(v, g + s)|ds
$$
$$
B(v,g):= \frac{1}{H} \left| \int_{g}^{g+H} [ \phi(v+h_1,w)-\phi(v,w) ] dw \right|, \qquad w=g+s.
$$
We first control the term $B$ as follows
$$
B(v,g)\leq \frac{1}{H}  \left( \left| \int_{g+H}^{+\infty} [ \phi(v+h_1,w)-\phi(v,w) ] dw \right|+  \left| \int_{g}^{+\infty}  [ \phi(v+h_1,w)-\phi(v,w) ] dw \right| \right).
$$
Applying inequality (\ref{estdvnum}) to each term, we obtain that
$$
\int_{\Omega_h}B(v,g) dvdg \leq 2C \frac{|h_1|}{H} \leq 2C|h_1|^{\frac{1}{2}} .
$$
Next, we control the term $A$. The Fubini Theorem and estimate (\ref{estdgnum}) give
$$ 
\int_{\Omega_h}A(v,g) dvdg \leq 2 C \frac{1}{H} \int_{0}^{H} sds  \leq CH \leq C |h_1|^{\frac{1}{2}}
$$
which concludes the proof of estimate (\ref{regest1}). 
\hfill $\square$

%%%%%%%%%%%%%%%%%%%%%%%%%%%%%%%
%-------------------------------------------------------------------------------------------
\subsection{Integrability with a singular weight}
%-------------------------------------------------------------------------------------------
%%%%%%%%%%%%%%%%%%%%%%%%%%%%%%%

Sobolev injections with the Besov regularity  (\ref{regest1})  imply (see  \cite{triebel, albert} for instance)
\beq
\big(- g_L v+g(V_E-v) \big) p(v,g) \in L^q \qquad \forall q, \quad 1 \leq  q < \frac{4}{3}.
\label{sobolev}
\eeq
The difficulty  in proving estimate (\ref{regest2}) is that the weight $J_v(v,g)=- g_L v+g(V_E-v) $ vanishes on a curve for $0\leq g \leq g_F$ and it remains  to prove that for some $G > g_F$ we have
\beq
\int_{0}^{G} \int_{0}^{V_F} p^q dvdg<+ \infty.
\label{eq:localint}
\eeq

A preliminary step in this direction is the following Lemma related to the hypoelliptic character  of equation (\ref{eq:iftd}). 

%----------------------------------------
\begin{lemma}\label{estlp}
For $0 \leq \al < 1$, we have 
\begin{equation}\label{estdgm}
 \int_{0}^{G} \int_{0}^{V_F}  \f{p(v,g)}{|J_v(v,g)|^\al} dg dv \leq C(\al, G).
\end{equation}
\end{lemma}
%-----------------------------------------

\proof 
Using that $p$ vanishes for $g=\infty$, we have
$$
p(v,g)=- 2 \int_g^\infty \nabla_g \sqrt{p(v,g')} \;  \sqrt{p(v,g')} dg' \leq 2\left[  \int_g^\infty \left(\nabla_g \sqrt{p(v,g')} \right)^2dg' \int_g^\infty  p(v,g') dg' \right]^{1/2} .
$$
Therefore
$$
 \int_{0}^{G}  \f{p(v,g)}{|J_v(v,g)|^\al} dg \leq  2\left[  \int_0^\infty \left(\nabla_g \sqrt{p(v,g')} \right)^2dg' \int_0^\infty  p(v,g') dg' \right]^{1/2}
\int_0^G \f{1}{|J_v(v,g)|^\al} dg.
$$
For the range of $\al$ under consideration, the last integral is bounded and thus, using the Cauchy-Schwarz inequality we find
$$\begin{array}{rl}
\dis \int_{0}^{V_F}  \int_{0}^{G}  \f{p(v,g)}{|J_v(v,g)|^\al} dg dv \leq & C(\al, G) \dis \int_{0}^{V_F} \left[  \int_0^\infty \left(\nabla_g \sqrt{p(v,g')} \right)^2dg' \int_0^\infty  p(v,g') dg' \right]^{1/2} dv
\\ \\
& \leq C(\al, G) \dis \int_{0}^{V_F} \int_0^\infty \left(\nabla_g \sqrt{p(v,g')} \right)^2dg' dv \;  \int_{0}^{V_F} \int_0^\infty  p(v,g') dg'dv.
\end{array}
$$
\hfill $\square$

%%%%%%%%%%%%%%%%%%%%%%%%%%%%%%%
%-------------------------------------------------------------------------------------------
\subsection{Proof of $L^q$ integrability}
%-------------------------------------------------------------------------------------------
%%%%%%%%%%%%%%%%%%%%%%%%%%%%%%%
 \bigskip

We are already reduced to proving \eqref{eq:localint} for some $G>g_F$. Take $1\leq q < 8/7$, we have 
$$\begin{array}{rl}
 \dis  \int_{0}^{G} \int_{0}^{V_F} p(v,g)^q dg dv & =  \int_{0}^{G} \int_{0}^{V_F} ( J_v^\al p)^{q/2} \left(\f{p }{|J_v|^\al} \right)^{q/2}dg dv
\\ \\
& \leq \left( \dis \int_{0}^{G} \int_{0}^{V_F} ( J_v^\al p)^{qr/2}dg dv \right)^{1/r} \left(  \dis\int_{0}^{G} \int_{0}^{V_F}  \left(\f{p }{|J_v|^\al} \right)^{qr'/2}dg dv \right)^{1/r'} 
\end{array} $$
after using H\"older's inequality.

We choose $qr'=2$, that is $qr=2(r-1)$, so that the last term is controled thanks to \eqref{estdgm}.

We claim that as long as $\f{qr}{2} < \f 4 3 $ the first term is also controled for $\al$ close enough to $1$ thanks to \eqref{sobolev}. This is because we can use again H\"older's inequality and write
$$
 \int_{0}^{G} \int_{0}^{V_F} ( J_v^\al p)^{qr/2}dg dv \leq  \left( \int_{0}^{G} \int_{0}^{V_F} ( J_v p)^{sqr/2}dg dv \right)^{1/s} 
  \left( \int_{0}^{G} \int_{0}^{V_F} \f{1}{J_v^{(1-\al) qrs'/2 } } dg dv \right)^{1/s'}. 
$$
We may choose $s>1$ such that we still have $sqr/2 < 4/3$ and thus use \eqref{sobolev} to bound the first term in the right hand side. On order to bound the second term, $qrs'/2$ is large but we may always choose $\al$ close enough to $1$ so that $(1-\al) qrs'/2<1$.
\\

It remains to compute the range of possible coefficients in view of these constraints; these are
$$
 qr=2(r-1), \qquad \f{qr}{2} < \f 4 3 .
$$
As announced, this gives  $1 \leq r < 7/3$ and $q < 8/7$. 
\hfill $\square$

%%%%%%%%%%%%%%%%%%%%%%%%%%%%%%%%%%%%%%%%
%-------------------------------------------------------------------------------------------------------------------
\subsection{Existence for the linear stationary equation}
%-------------------------------------------------------------------------------------------------------------------
%%%%%%%%%%%%%%%%%%%%%%%%%%%%%%%%%%%%%%%%

Here we sketch a path to prove existence but we do not go to the full details which would need another paper. 

With the bounds at hand, there are several ways to obtain existence of a steady state. A method is to reduce the problem to a finite dimensional linear system with positivity and apply the Perron-Frobenius theorem. This is equivalent to write a stable (positivity preserving) numerical discretization for the equation \eqref{eq:iftd}.

To do so, it is usual to symmetrize the diffusive part of the equation, use the unknown 
$$q (v,g) = e^{\f{(g-g_{\rm in})^2}{2a}}p(v,g)$$
  and take a bounded domain to obtain the following equation on $q$
\begin{equation}
\label{eq:ifsym}
   \frac{\partial}{\partial v} \left[J_v e^{-\f{(g-g_{\rm in})^2}{2a}}q(v,g)\right] 
-\frac{a}{\sigma_{E}}  \frac{\partial}{\partial g} \left[ e^{-\f{(g-g_{\rm in})^2}{2a}} \frac{\partial}{\partial g}  q(v,g)\right]   =  0, \qquad 0\leq v \leq V_F, \; 0 \leq g \leq R, 
\end{equation}
with the normalization $\int_0^{V_F} \int_0^R e^{-\f{(g-g_{\rm in})^2}{2a}}q(v,g) dv dg =1$ and the zero flux boundary conditions
$$\bepa
q(0,g) = q(V_F,g)  = 0, \qquad &\text{  for  } \;  g  \leq g_F,
\\ \\
g V_E q(0,g) = [- g_L V_F +g(V_E-V_F) ] q(V_F,g) , \qquad &\text{ for } \; g > g_F ,
\\ \\
 \frac{\partial}{\partial g}  q(v,g=0) =  \frac{\partial}{\partial g}  q(v,g=R) =0, \qquad & \text{ for } \;0 \leq v \leq V_F,
\eepa$$
A first step is to obtain the same estimates as above for this problem and show that, as $R \to \infty$ we obtain the solution of  \eqref{eq:iftd}.

The second step is to build a finite dimensional approximation of \eqref{eq:iftd}. This can be achieved by a standard finite volume scheme (see for instance \cite{bouchut_book, CCTao}). It uses a discretization step $h$ (which we take the same in $v$ and $g$ to simplify), the  discrete grid $v_{i+ \f 12} =  ({i+ \f 12} ) hV_F$, $0 \leq i \leq I$ with $Ih=1$ and  $g_{j+ \f 12} = ({j+ \f 12}) h R $, $0 \leq j \leq I$. The solution $q$ is approximated by 
$$
h^2 V_F R \; q_{i,j} \approx  \int_{v_{i-\f 12}}^{v_{i+ \f 12}}  \int_{g_{i-\f 12}}^{g_{i+ \f 12}} q(v,g) dv dg.
$$
The finite dimensional problem is then written using an upwind discretization of the first order derivatives and a centered scheme for the diffusion term, that is for $1\leq i,\; j \leq I$, 
\begin{equation}
\label{eq:ifnum}
\f{1}{hV_F} \left[ {\mathcal J}_{i+\f 12,j} q_{i+\f 12,j} -    {\mathcal J}_{i-\f 12,j} q_{i-\f 12,j} \right]
- \frac{a}{  h^2R^2 \sigma_{E}} \left[  {\mathcal M}_{j+\f 12} ( q_{i,j+1} - q_{i,j} )  -  {\mathcal M}_{j -\f 12} ( q_{i,j} - q_{i,j-1 })  \right]   =  0,
\end{equation}
with 
$$
 {\mathcal M}_{j+\f 12}= e^{- (g_{j+\f 12}  - g_{\rm in})^2/(2a)} , \qquad  \qquad {\mathcal J}_{i+\f 12,j} = \left[ -g_L v_{i+\f 12} +g(V_E- v_{i+\f 12}) \right]e^{- (g_{j}  - g_{\rm in})^2/(2a)} ,
$$
$$
q_{i+\f 12,j}  = q_{i,j} \text{ if }Ê  {\mathcal J}_{i+\f 12,j}\geq 0, \qquad = q_{i+1,j} \text{ if }Ê  {\mathcal J}_{i+\f 12,j}\leq 0.
$$
The advantage of this approach is that the boundary conditions come naturally as $q_{i,0 }= q_{i,1 }$ and $q_{i,I+1 }= q_{i,I }$ and define the endpoint terms on the diffusion approximation. Also it is easy to check that for the $j$ such that $g_j \leq g_F$ one should take (and needs to define) $ q_{\f 12,j} =0= q_{I+\f 12,j}$. And for the $j$ such that $g_j \geq g_F$ one takes 
$  {\mathcal J}_{\f 12,j} q_{\f 12,j}  = {\mathcal J}_{I +\f 12,j} q_{I+\f 12,j} $ because $q_{\f 12,j} $ is not defined but $q_{I+\f 12,j}=q_{I,j}$.

One readily checks that this problem corresponds to a matrix with positive diagonal and nonpositive terms out of the diagonal. The diagonal is dominant but not strictly and $0$ is the first eigenvalue with a positive vector in its kernel. Again the lengthy calculation is to redo at the discrete level the above estimates so as to prove bounds better than $\ell^1$ bounds and to pass to the limit as $h \to 0$.

%%%%%%%%%%%%%%%%%%%%%%%%%%%%%%%%%%%%%%%%
%-------------------------------------------------------------------------------------------------------------------
\subsection{Uniqueness of the linear stationary equation}
%-------------------------------------------------------------------------------------------------------------------
%%%%%%%%%%%%%%%%%%%%%%%%%%%%%%%%%%%%%%%%

%-------------------------------------------------
 \begin{lemma}[Uniqueness]
 A   nonnegative solution of equation \eqref{eq:iftd} with the boundary conditions  \eqref{eq:BC1}, \eqref{eq:BC2}, \eqref{eq:NF} is unique.
\label{lem:unicite}
\end{lemma}
%-------------------------------------------------
%
{\bf Proof of Lemma~\ref{lem:unicite}.}
Let $p_1$ and $p_2$ be two nonnegative solutions of equation  \eqref{eq:iftd}. We set 
$$
 p_m:=\f{p_1 + p_2}{2} ,
$$
another nonnegative solution for which $ P:= \frac{|p_1-p_2|^2}{p_m}$ is well defined. 

Then, the function $P$ satisfies the relative entropy equality (see \cite{Pe} pp 166--167)
$$ 
\frac{p_m}{\sigma_E}  \left(\partial_g\frac{|p_1-p_2|}{p_m} \right)^2 + \partial_v [(-g_Lv+g(V_E-v))P] +  \frac{1}{\sigma_E} \left(\partial_g [(-g +g_{\rm in}) P] -  a  \partial_g^2P \right) = 0.
$$
As $p_1$ and $p_2$ satisfy  the boundary condition \eqref{eq:NF},  one readily checks that the no-flux condition also holds 
$$
(-g +g_{\rm in})P - a  \partial_gP =0 \qquad \hbox{ at } \; g=0, \quad g=+\infty.
$$
In the same way, from (\ref{eq:BC1}), (\ref{eq:BC2}), we observe that $P$ also satisfies them also
$$
J_v(0,g)P(0,g)=J_v(V_F,g)P(V_F,g), \qquad \forall g \geq 0.
$$

We may now integrate the equation on $P$ and use the boundary conditions (\ref{eq:BC1}), (\ref{eq:BC2}), (\ref{eq:NF}) on $P$ to obtain
$$
\int_{0}^{+\infty} \int_{0}^{V_F}p_m  \left(\partial_g\frac{|p_1-p_2|}{p_m} \right)^2 dvdg=0,
$$  
from which we deduce that $\frac{p_1-p_2}{p_m}=: Q(v)$ is independent of $g$. Using that $p_1-p_2$ and $P=Q(v)(p_1-p_2)$ satisfy 
$$ 
 \partial_v [(-g_Lv+g(V_E-v))P] +  \frac{1}{\sigma_E} \left(\partial_g [(-g +g_{\rm in}) P] -
a  \partial_g^2P \right) = 0, 
$$
it follows easily that $Q$ is a constant and Lemma~\ref{lem:unicite} is proved. 
\hfill $\square$   

%------------------------------------------------------------------------ 
\subsection{Long time convergence to the steady state}
%------------------------------------------------------------------------

 As a consequence of our analysis of the steady state equation, we prove some long time behavior for the solution of \eqref{eq:iftdp} in the linear case $S_E=0$. More precisely, we prove that the integral of the solution with respect to the variable $v$, converges  exponentially converges to the stationary state. 
 Indeed, in the linear case, after integrating in $v$,  the stationary state is then an explicit Maxwellian (see section  \ref{sec:stst}) and we find that  $\int_{0}^{V_F} p(v,g,t)dv$
 is solution of a  classical one dimensional Fokker-Planck equation. Hence,  we are able to use a generalized entropy method associated to a Poincar{\'e} inequality in the spirit of the method describe in \cite{Pe, CCP}. These arguments are also used by  \cite{calvezmeunier, lepoutremeunier} for the study of Fokker-Planck equations where  blow-up is proved and where generalized entropy method with Poincar{\'e} inequality is used. The following Theorem holds
 %----------------------------------------------
 \begin{thm}\label{thm:linear}
 For a global solution of equation \eqref{eq:iftdp} with $S_E=0$,  let 
 $$M(g)= Z(g_{\rm in})^{-1} e^{- \frac{(g-g_{\rm in})^2}{2a}} \quad \hbox{ and  } \; \varphi(g,t):=\int_{0}^{V_F} p(v,g,t) dv.$$ 
 Then, there exists $\eta>0$ such that
 $$ \int_{0}^{+\infty} M(g) \left( \frac{\varphi}{M}\right)^2 (g,t)dg \leq  e^{- \eta t}  \int_{0}^{+\infty} M(g) \left( \frac{\varphi}{M}\right)^2 (g,0)dg.$$ 
\end{thm}
%---------------------------------------------------- 
{\bf Proof of Theorem~\ref{thm:linear}.}
Integrating equation  \eqref{eq:iftdp}  between $0$ and $V_F$, we find that $\varphi$ is solution of the equation
$$\partial_t \varphi  +\frac{1}{\sigma_{E}} \left(\frac{\partial}{\partial g} \left[ ( g_{\rm in} - g ) \varphi \right]  
- a \frac{\p^2}{\p g^2} \varphi \right) =  0.$$
After classical computations (see  for example book \cite{Pe} pp166-167),  we find that
$$\frac{d}{dt}  \int_{0}^{+\infty} M(g) \left( \frac{\varphi}{M}\right)^2 dg\leq -  \frac{a }{\sigma_E} \int_{0}^{+\infty} M(g) \left(\partial_g \left( \frac{\varphi}{M}\right)\right)^2 dg.$$
As $M(g)$ is a Maxwellian with
$$\int_{0}^{+\infty}M(g) dg=1,$$
we can apply   Poincar{\'e} inequality (see \cite{Pe} page 167, \cite{CCP} and references therein) to obtain that there exists a constant $C>0$ such that
$$   \int_{0}^{+\infty} M(g) \left( \frac{\varphi}{M}\right)^2 dg \leq C \int_{0}^{+\infty} M(g) \left(\partial_g \left( \frac{\varphi}{M}\right)\right)^2 dg.$$
So, there exists a constant $\eta>0$ such that
$$\frac{d}{dt}  \int_{0}^{+\infty} M(g) \left( \frac{\varphi}{M}\right)^2 dg\leq -\eta   \int_{0}^{+\infty} M(g) \left( \frac{\varphi}{M}\right)^2 dg$$
and Theorem~\ref{thm:linear} follows thanks to the Gronwall lemma. 
\hfill $\square$

%%%%%%%%%%%%%%%%%%%%%%%%%%%%%%%%%%%%
%--------------------------------------------------------------------------------------------------------
\section{Steady states for the nonlinear equation}
\label{sec:ststnl}
%--------------------------------------------------------------------------------------------------------
%%%%%%%%%%%%%%%%%%%%%%%%%%%%%%%%%%%%

Our analysis of the linear equation has immediate consequences on the nonlinear case. 
We still assume that the external input rate $\nu>0$  is  constant. Notice that, when deriving of the kinetic equation,   this  corresponds to a network in which the input spikes follow a Poisson process.

We recall the $g$-dependent firing rate, the total firing rate, the input current and the noise are defined as 
\beq
N(g):= [- g_L V_F +g(V_E-V_F) ] p(V_F,g) \geq 0, \qquad   {\mathcal N} :=  \int_0^\infty  N(g) dg,
\label{recall:nl}
\eeq
\beq
g_{\rm in} = f_E \nu+ S_E {\mathcal N}, \qquad \qquad a= \frac{1}{2 \sigma_{E}} \left( f_E^{2} \nu +   \f{ S_{E}^2}{N_E} \mathcal{N}\right).
\label{recall:nl2}
\eeq

The aim of  this section is to prove the following Theorem 
%---------------------------------------------------------
\begin{thm}\label{thm:fixedpoint}
When  the strength of interconnections $S_E$ is  such that
$$
\frac{V_E}{V_F} S_E <1  \qquad \text{(weakly connected network)  }
$$
there exists at  least one  solution to \eqref{eq:iftd}, \eqref{recall:nl}, \eqref{recall:nl2} with the same regularity as in Theorem~\ref{thm:statmain}.

When the conditions  hold
$$
\frac{V_E-V_F}{V_F} S_E >1 \quad  \hbox{ and } \quad  (V_E-V_F)  f_E \nu > V_F^2  \qquad \text{(strong connection,  strong noise),  }
$$
equation \eqref{eq:iftd}, \eqref{recall:nl}, \eqref{recall:nl2}  does not have solutions.
\end{thm}
%---------------------------------------------------------
% 
{\bf Proof of Theorem~\ref{thm:fixedpoint}.} 
We define the following function 
\begin{equation}\label{def:eqpsi}
{\mathcal N} \geq 0 \mapsto  \Psi( {\mathcal N} ):= V_E\int_0^{+\infty}g \; p(0,g) dg= \int_0^\infty N(g) dg  \geq 0
\end{equation}
where $ p(v,g) $ is the solution of  equation \eqref{eq:iftd} with  \eqref{recall:nl2} for this given ${\mathcal N}$. A solution to the nonlinear problem is a fixed point of $\Psi$.
\\

We first prove some properties on $\Psi$ which are a first step toward our goal
%------------------------------------------- 
\begin{lemma}\label{lem:psifp}
Let $\Psi$ defined as in \eqref{def:eqpsi}. Then the following properties hold
\\[5pt]
(i) \ $\Psi(0)$ is positive, 
\\[5pt]
(ii)  $\Psi$ is continuous,
\\[5pt]
(iii) \vspace{-10pt}
\begin{equation}\label{eqbounpsi}
-V_F+ (V_E-V_F) \left( g_{\rm in} + \sqrt{ \frac{a }{2\pi}} e^{- \frac{ g_{\rm in}^{2}}{2a} } \right)\leq V_F \Psi( {\mathcal N} ) \leq  V_E  \left( g_{\rm in} + \sqrt{ \frac{2 a }{\pi}} e^{- \frac{ g_{\rm in}^{2}}{2a} }\right).
\end{equation}
\end{lemma}
%-------------------------------------------

\noindent {\bf Proof of Lemma~\ref{lem:psifp}.}
For the item (i), we prove that for $\mathcal{N}=0$, the solution $p$ does not vanish on  the boundaries $v=0$ or $v=V_F$. Otherwise it would satisfy a $0$ Dirichlet condition for which the only solution vanishes everywhere. This can be seen thanks to the dual equation which has a super-solution. A direct way to see this is to integrate by parts \eqref{eq:iftd} against a weight $\Phi(v,g)$; if $p$ vanishes at $v=0$ and $v=V_F$ we have
$$\bea
\int_0^\infty \int_0^{V_F} p(v,g) \big[ -\p_v \Phi (v,g) J_v(v,g)& + \left(\f{g- f_E \nu}{\sigma_E}\p_g \Phi (v, g) -  \frac{a}{\sigma_E} \p^2_{gg} \Phi (v,g) \big] \right) dv dg 
\\ \\
&=  \frac{a}{\sigma_E} \int_0^{V_F} \p_g \Phi (v, 0) p(v,0) dv.
\eea $$
We choose the weight $\Phi = e^{\lb g} e^{\mu v}$ with $\lb >0$, $\mu >0$ and find 
$$
\int_0^\infty \int_0^{V_F} p(v,g)  \Phi (v,g) \big[ \mu (g_Lv-g(V_E-v)) +  \frac{1}{\sigma_E}(\lb(g-\nu) - a\lb^2 )\big] dv dg  \geq 0. 
$$
Therefore 
$$
\int_0^\infty \int_0^{V_F} p(v,g)  \Phi (v,g) \big[ \mu (g_LV_F-g(V_E-V_F)) +  \frac{1}{\sigma_E}(\lb(g-\nu) - a\lb^2 )\big] dv dg  \geq 0. 
$$
We now take $\mu (V_E-V_F)= \frac{1}{\sigma_E} \lb$ and arrive to 
$$
\int_0^\infty \int_0^{V_F} p(v,g)  \Phi (v,g)  \lb \big[ \f{g_LV_F}{V_E-V_F}  + \frac{1}{\sigma_E}(- \nu -a \lb) \big] dv dg  \geq 0.
$$
For $\lb $ large enough, this means that $p=0$ which contradicts the probability normalization. 
\\

\noindent Item (ii)  is a consequence of uniqueness for  solutions of \eqref{eq:iftd}.  Indeed, let  $({\mathcal N}_k)_{k \in \N}$ a sequence which converges to  ${\mathcal N}$ and let $(p_k)_{k \in \N}$ be the unique solution of equation \eqref{eq:iftd} with this input  ${\mathcal N}_k$.  Let $\widetilde{p}_k(g)$ be the Maxwellian associated to $p_k$ given, according to \eqref{eq:intv}, by 
$$
 \widetilde{p}_k (g):= \int_{0}^{V_F} p_k(v,g) dv = Z_k^{-1} e^{-\frac{(g-g_{\rm in})^2}{2a_k}} .
$$
Because $Z_k$ and $a_k$ converge,  $\widetilde{p}_k$  converges strongly  in $L^1(0,\infty)$ to the maxwellian  $\widetilde{p}$ and thus
\begin{equation}\label{est2} 
\int_0^{+\infty} \widetilde{p}(g)  dg=1.
\end{equation}
On the other hand, we know that
$$
\int_0^{+\infty}  \int_{0}^{V_F} p_k(v,g) dvdg=1,
$$
and we deduce with estimate (\ref{regest2}) that, modulo a subsequence, $p_k$ converges weakly in $L^q$ to a function $p$ such that
$$
\int_0^{+\infty}  \int_{0}^{V_F}p(v,g) dv dg=1.
$$
Passing to the limit in equation \eqref{eq:iftd}, we obtain that $p$  is the normalized solution of 
 \eqref{eq:iftd}  associated to  ${\mathcal N}$ which finishes the proof of continuity of $\Psi$. 
\\ 
 
\noindent For (iii), we write
$$
V_F \Psi( {\mathcal N} ) = \int_0^{V_F}  \int_0^\infty [-g_L v+ g(V_E-v)] p(v,g)dgdv
$$ 
and so
$$  
-V_F+ (V_E-V_F)  \int_0^{V_F}  \int_0^\infty gp(v,g)dgdv \leq V_F \Psi( {\mathcal N} ) \leq  V_E \int_0^{V_F}  \int_0^\infty gp(v,g)dgdv.
$$
We conclude estimate (\ref{eqbounpsi}) with inequality  \eqref{estgin} and the proof of Lemma~\ref{lem:psifp} is complete.   \hfill $\square$
 
 \vspace{0,5cm}
 
We now have the material  to conclude the proof of Theorem~\ref{thm:fixedpoint}. For weak interconnections, as $\Psi(0)>0$, and $\Psi$ continuous (see Lemma~\ref{lem:psifp}), to prove existence of at least one steady state, it is enough to know that  there exists $x \in \R^+$ such that
$$
\Psi(x)-x<0.
$$
For this, we use the right hand side of estimate \eqref{eqbounpsi}, the explicit formulas on $a$ and $g_{\rm in}$ and  estimate (\ref{eqbounpsi}) on $\Psi$ to 
 obtain that there exists two constants $C_1>0$ and $C_2>0$  such that 
$$
\Psi(x)-x \leq x \left( \frac{V_E}{V_F}  S_E -1 \right)+ C_1+ C_2 \sqrt{x}
$$
and so,  under the condition $ \frac{V_E}{V_F}  S_E <1$, we conclude that 
$$
 \lim_{x \to +\infty} \Psi(x)-x =-\infty.
$$
This proves the first point of Theorem~ \ref{thm:fixedpoint}.
\\

For strong interconnections, using now the left hand side of estimate \eqref{eqbounpsi}, we obtain that 
$$
\Psi(x)-x \geq -V_F+   \frac{V_E-V_F}{V_F}  \left( f_E \nu+ S_E x+  \sqrt{\f{a}{2 \pi} }e^{- \frac{ g_{\rm in}^{2}}{2a} }\right) -x
$$
which, with our strong connection and strong noise assumptions, implies that
$$
\Psi(x)-x>0
$$
which means that there is no stationary solution and the proof of Theorem~  \ref{thm:fixedpoint} is complete.
\hfill $\square$
 
%-------------------------------------------------------------------------------------------------------------------------------------------
\section{Evolution equation: estimates on moments and on the firing rate}
\label{sec:evol_moments}
%-------------------------------------------------------------------------------------------------------------------------------------------
 
We begin our analysis of a priori estimates of the evolution equation \eqref{eq:iftd}  by showing that the solution propagates several types of moments in $g$ in $L^1$.  Moreover, we obtain a priori $L^{1}_{loc}$ bounds on the total firing rate $\mathcal{N}(t)$ which on the one hand  implies that, for high interconnections, the firing rate cannot be bounded, for initial datas    localized in high  $g$ and on the other hand are the first steps towards $L^q$ estimates, $q>1$ which are performed in the next section. These estimates rely on the combination of moments as follows.

%------------------------------------------------------------------------  
\subsection{Equations on the moments}
\label{sec:moments} 
%------------------------------------------------------------------------

For $k \in \N$ and we use the notations
$$
 \psi(t):=  \int_{0}^{V_F} v p(v,g,t)dvdg, \qquad h_k(t):=  \int_{0}^{V_F} \int_{0}^{+\infty} g^kp(v,g,t)dvdg, \qquad f(t):= \int_{0}^{V_F}p(v,0,t)dv.
$$
We recall that, $p(t)$ being a probability density, we have
$$
h_1(t)^2 \leq h_2(t).
$$

 We assume that
\begin{equation}\label{initial_k}
\int_{0}^{V_F}\int_{0}^{+\infty}(1+g)^{k}p(v,g,0)dvdg<+\infty.
\eeq
Then, one readily checks the following differential relations
\begin{equation}\label{eqpsi}
 \frac{d}{dt} \psi(t)=  - g_L \psi (t)  - \int_{0}^{V_F} \int_{0}^{+\infty} gv \;  p(v,g,t) dvdg +V_E h_1(t) - V_F \mathcal{N}(t),
\end{equation}
\begin{equation}\label{eqh1}
 \frac{d}{dt} h_1(t) =   \frac{1}{\sigma_E} [ -  h_1(t) + g_{\rm in}(t) +a(t)f(t) ],
\end{equation}
\begin{equation}\label{eqh2}
\frac{d}{dt}h_2(t) = \frac{1}{\sigma_E} [ -2 h_2(t)+ 2 g_{\rm in}h_{1}(t)+2a(t) ],
\end{equation}
and more generally for $k \geq 2$, 
\begin{equation}\label{eqhk}
\frac{d}{dt}h_k(t) = \frac{1}{\sigma_E} [ -k h_k(t) + k g_{\rm in}h_{k-1}(t)+k(k-1)a(t) h_{k-2}(t) ]
\end{equation}

These are obtained integrating equation (\ref{eq:iftdp}) by parts after multiplying it respectively by the weights   $v$, $g$, $g^2$  and  $g^k$ for $k \geq 2$.
\\

 In order to manipulate these relations, the difficulty is that neither $f(t)$ nor $ \mathcal{N}(t)$ are under control.  Our  goal is first to explain how we can go around it.

%%%%%%%%%%%%%%%%%%%%%%%%%%%%%%%
%-------------------------------------------------------------------------------------------
\subsection{Upper bounds on the moments $h_k$}
%-------------------------------------------------------------------------------------------
%%%%%%%%%%%%%%%%%%%%%%%%%%%%%%%

 We now derive the a priori bounds which show that whatever the strength of interconnections, if the initial data has  $k \geq 2$ moments, then so does the solution for all time. Moreover, if the strength  of interconnections $S_E$ is small enough, then we have a uniform bound $h_1 \in L^{\infty}(\R^+).$ 

It is convenient to define the  two numbers $\lb_E >0$ and $\omega_E \in \R $ by
\beq
\lb_E V_F =\frac{S_E}{\sigma_E} + \frac{ S_E^2 }{2 N_E \sigma_{E}}, \qquad \quad \sigma_E \omega_E= \f{V_E}{V_F}  \left[ S_E +\frac{ S_E^2 }{2 N_E }\right] - 1 .
\label{def:param}
\eeq
Then, we can establish the controls 
%--------------------------------------------------- 
\begin{thm}\label{thm_upper}  
Let $k \geq 2$ and assume that the initial data satisfies \eqref{initial_k}. 
Then, there is a constant $C$ such that the following a priori estimates hold
\begin{equation}\label{estthmh1upper}
 h_1(t) =  \int_{0}^{V_F} \int_{0}^{+\infty} gp(v,g,t)dvdg  \leq C  \max(1,e^{ \omega_E t}),
\end{equation}
\begin{equation}\label{estthmfupper}
  \int_0^t  \mathcal{N}(s)ds\leq C(1+t)  \max(1,e^{ \omega_E t})   .
\end{equation}  
Moreover, for all $T>0$ there is a constant $C(T)$ such that
\begin{equation}\label{esthk}
\sup_{t \in [0,T]}   \int_{0}^{V_F} \int_{0}^{+\infty} g^kp(v,g,t)dvdg\leq C(T) , \qquad \int_{0}^{T} \int_{0}^{+\infty} (1+g)^{k-1} N(g,t) dg\leq C(T).
\end{equation}
\end{thm}
%---------------------------------------------------

The regime $\omega_E<0$ corresponds again to weak interconnections, i.e. $S_E$ small, however with a different definition than for steady states in section \ref{sec:ststnl} which is sharper. Then, \eqref{estthmh1upper} and \eqref{estthmfupper} give uniform controls in $L^\infty$.  
\\

\noindent {\bf Proof of Theorem~\ref{thm_upper}.} 
The  estimates  (\ref{estthmh1upper}) and (\ref{estthmfupper}) come together by a combinations of the relations in section \ref{sec:moments}.

Firstly, we multiply equation (\ref{eqh1}) by $h_1$ and subtract it to  equation (\ref{eqh2}). We find that 
$$\displaystyle \bea
\frac{d}{dt}Ê\f{h_2-h_{1}^2}{2} &=\frac{1}{\sigma_E}\left[- (h_2-h_{1}^2)+ a(t) - a(t) f(t) h_1(t) \right]
\\ \\
& \leq - \frac{2}{\sigma_E} \f{h_2-h_{1}^2}{2} - \f{a(t) f(t) h_1(t)}{\sigma_E} + \f{f_E^2 \nu_M} {2\sigma_E}+ \f{S_E^2}{2\sigma_E N_E}  \mathcal{N}(t).
\eea$$
The last line is a consequence of formula (\ref{defa}) on  $a$ and of assumption (\ref{asnu}).

Secondly, we define the function $G(h) =(h-1)_+$. Using respectively that $G'(h) \leq h$, $hG'(h) \geq G(h)$ and $G'(h) \leq 1$, we have from \eqref{eqh1}
$$\bea
\frac{d}{dt}G(h_1) &=   \frac{1}{\sigma_E} \left[ G'(h_1) a(t)f -G'(h_1)h_1+ g_{\rm in}(t)G'(h_1)\right)
\\ \\
& \leq  \frac{1}{\sigma_E} \left[ h_1(t) a(t)f(t) -G(h_1)+ g_{\rm in}(t) \right]
\\ \\
& \leq  \frac{1}{\sigma_E} \left[ h_1(t) a(t)f(t) -G(h_1)+ f_2 \nu_M+S_E \mathcal{N}(t) \right].
\eea$$
The last line is a consequence of the definition of $g_{\rm in}$ in (\ref{defdrift}) and of assumption (\ref{asnu}).
\\

Thirdly, equation (\ref{eqpsi}) gives 
$$ 
\frac{d}{dt} \psi(t) \leq V_E h_1- V_F \mathcal{N}(t) \leq V_E G(h_1) - V_F \mathcal{N}(t)+ V_E  .
$$

We can form the combination of three nonnegative quantities 
$$
k(t):= \f{h_2-h_{1}^2}{2} + G(h_1)+ \lambda_E \psi 
$$ 
which satisfies
\begin{equation}\label{eq:k}
\frac{d}{dt} k(t) \leq   - \frac{2}{\sigma_E} \f{h_2-h_{1}^2}{2} + \omega_E G(h_1) + C.
\end{equation}
We continue our proof  differently in the cases when $\omega_E \geq 0$ or 
$\omega_E<0$.
\\

\noindent  {\it $\bullet$ Case when $\omega_E\geq 0$.}
Because $0 \leq \psi \leq V_F$, using \eqref{eq:k} there is a constant  $C$ such that
\begin{equation}\label{eq:kp}
\frac{d}{dt} k(t) \leq  \omega_E k(t) +C .
\end{equation}
We can use the Gronwall lemma and obtain
$$
k(t) \leq C e^{\omega_E t} 
$$
which gives (\ref{estthmh1upper}). For  (\ref{estthmfupper}), we integrate equation \eqref{eqpsi}  on $0\leq \psi(t) \leq V_F$ to obtain that (the case $\omega_E=0$ is treated to the expense of a factor $1+t$)
$$
\int_{0}^{t} \mathcal{N}(s) ds \leq  \psi(0)+ V_E \int_0^t h_1(s) ds  \leq C(1+t)  e^{\omega_E t}.
$$

\noindent $\bullet$ {\it Case where $\omega_E<0$.}
Then,  $\omega_E \geq -\frac{2}{\sigma_E}$. Hence, we can use again equation  \eqref{eq:k} to conclude that  \eqref{eq:kp} still holds true and thus $k(t) \leq C$. This concludes the proof of the bounds  (\ref{estthmh1upper}) and (\ref{estthmfupper}). 
\\

We now come to (\ref{esthk}). The first estimate on $h_k$ is a consequence of equation (\ref{eqhk}) when iterating on $k$, after integration in time and using estimates (\ref{estthmh1upper}) and (\ref{estthmfupper}).

To prove the second bound, we multiply equation \eqref{eq:iftdp} by $vg^{k}$ and find that
$$\bea
\frac{d}{dt} \int_{0}^{V_F}\int_{0}^{+\infty} vg^k p & dvdg  = -V_F \int_{0}^{+\infty} N(g,t)g^{k} dg + \int_{0}^{V_F}\int_{0}^{+\infty} g^{k} (-g_Lv+g(V_E-v)) p(v,g)dvdg
\\ \\
& + k \int_{0}^{V_F}\int_{0}^{+\infty}  g^{k-1} (-g+g_{\rm in}) p(v,g,t) dvdg + a k(k-1)\int_{0}^{V_F}\int_{0}^{+\infty} \int_{0}^{V_F} g^{k-2} p(v,g,t) dvdg.
\eea$$
Integrating the above equation, estimate  (\ref{estthmfupper}) and because we have already proved that $h_{k+1}$ is locally bounded,  we deduce that for all $T>0$
$$
\int_{0}^{T} \int_{0}^{+\infty} (1+g)^{k-1} N(g,t) dg \leq C(T)
$$
which ends the proof of Theorem~\ref{thm_upper}.
 \hfill $\square$

%----------------------------------------------------------------------------------------------------------
\subsection{Exponential growth on the moments for strong interconnections}
%----------------------------------------------------------------------------------------------------------

For high interconnections we can prove that the previous results are somehow sharp.   Indeed,  the first moment of the solution as well as  the total firing rate $\mathcal{N}(t)$ grow exponentially. This implies in particular that we cannot have periodic dynamic on the firing rate and that necessary the firing rate is not in $L^{\infty}(\R^+)$. The following theorem gives a precise condition on the coefficients
%
%----------------------------------------------------------------------------------------------------
\begin{thm}\label{thm:lowerbound}
Assume that the condition
$$
\zeta:=  \frac{1}{\sigma_E} \left(\frac{S_E(V_E-V_F)}{V_F}-1\right)>0
$$
 holds.
If  $h_1(0)$ is big enough, there exists a constant $C$ such that the following a priori estimates hold
$$ 
h_1(t) \geq Ce^{\zeta t}
$$
and for $t$ big enough
$$
\int_{0}^{t} \mathcal{N}(s)ds \geq Ce^{\zeta t}.
$$
\end{thm}
%-------------------------------------------------------------------------------------------------------

\noindent {\bf Proof of Theorem~\ref{thm:lowerbound}.}
From equations (\ref{eqpsi}) on $\psi$ and $(\ref{eqh1})$ on $h_1$ we infer  that
$$\frac{d}{dt} \psi \geq  (V_E- V_F) h_1 -g_L \psi -V_F \mathcal{N}(t)
\qquad 
\frac{d}{dt} h_1 \geq \frac{1}{\sigma_E} \left(   g_{\rm in}(t) -h_1   \right).
$$
To eliminate $ \mathcal{N}(t)$, we multiply equation on $\psi$ by $b=\frac{S_E}{V_F \sigma_E}$, use formula (\ref{defdrift}) for $g_{\rm in}$ and assumption (\ref{asnu});  we find that there is a constant $C>0$ such that
$$
\frac{d}{dt} (b\psi+h_1) \geq \left[ b(V_E-V_F)- \frac{1}{\sigma_E} \right] (b\psi+h_1)  - C= \zeta (b\psi +h_1)  -C .
$$
We obtain that 
$$
(b\psi+h_1)(t)\geq \left[ (b\psi+h_1)(0) - \frac{C}{\zeta} \right]e^{\zeta t} +\frac{C}{\zeta}.
$$
Since $\psi$ is bounded, we obtain the lower bound on $h_1$ just taking an initial data such that
$$
(b\psi+h_1)(0) >\frac{C}{\zeta} .
$$

To conclude the estimate on $\mathcal{N}(t)$, we come back to the equation on $\psi$ and use  this lower bound on $h_1$.
\hfill $\square$

%---------------------------------------------------------------------------------------------------------
\section{Integrability for the evolution equation and bounds on the firing rate}
\label{sec:reg}
%-------------------------------------------------------------------------------------------------------

Our next task is to obtain a priori estimates where higher integrability is propagated by the evolution equation. As a consequence we obtain that the firing rate $\mathcal N$ belongs locally to Lebesgue $L^q$ spaces for some $q>1$ thus discarding a blow-up scenario as in \cite{CCP}.

%------------------------------------------------------------------------
\subsection{Entropy inequality}
%------------------------------------------------------------------------

Using the $L^1$ control of  $g_{\rm in}$ that has been proved, we can obtain an entropy bound for the solution to the evolution equation
%--------------------------------------------
\begin{thm}\label{thm:regularity}
Assume that the initial data is such that
$$
 \int_{0}^{V_F} \int_{0}^{+\infty}\left[ | \ln(p(v,g,0)) | + g^2 +1\right]     p(v,g,0)  dvdg < \infty . 
$$
Then, for all $T>0$, the solution of equation   \eqref{eq:iftdp} satisfies the following a priori estimates
$$\bea
  \sup_{t \in [0,T]}\int_{0}^{V_F} \int_{0}^{+\infty} | \ln(p)|p(v,g,t) dvdg &+ \int_{0}^{T}\int_{ g_{F}} ^\infty N(g,s) \ln \f{g V_E }{ gV_E-gV_F- g_L V_F } dg ds 
\\ \\
 & + \f 1 {\sg_E}  \int_{0}^{T} \int_{0}^{V_F}\int_{0}^{+\infty}  a(s) \f{|\nabla_g p|^2}{p} dv dg  ds\leq C(T).
\eea $$
\end{thm}
%
%--------------------------------------------
%
{\bf Proof of Theorem~\ref{thm:regularity}.} 
We multiply equation   \eqref{eq:iftdp} by $\ln(p)$  and integrate to find that 
$$\bea
\frac{d}{dt} \int_{0}^{V_F} \int_{0}^{+\infty} & \ln(p) p(v,g,t) dvdg= -  \int_0^{V_F} \int_0^\infty  \left[- \big( g_L +g \big) p(v,g,t)  + \frac{1}{\sigma_E}a \f{|\nabla_g p|^2}{p}\right] 
\\ \\
&-  \frac{1}{\sigma_E}\int_0^\infty \big( -g_L v  +g(V_E -v )\big)  p) \ln p  \Big|_0^{V_F} \; dg
+  \int_0^\infty  \big( -g_L v  +g(V_E -v )\big)  p \Big|_0^{V_F} \; dg
\\ \\
& + \frac{1}{\sigma_E} \int_{0}^{+\infty}\int_0^{V_F}   \big( g_{\rm in} - g \big)  \partial_gp \; dgdv.
 \eea $$
Therefore, applying boundary conditions (\ref{eq:BC1}) and (\ref{eq:BC2}), we deduce that
$$ \bea
\frac{d}{dt} \int_{0}^{V_F} \int_{0}^{+\infty} & \ln(p) p(v,g,t) dvdg + \int_{ g_{F}} ^\infty N(g,t) \ln \f{g V_E }{ gV_E-gV_F- g_L V_F } dg  +  \int_{0}^{V_F} \int_{0}^{+\infty}  \frac{a}{\sigma_E} \f{|\nabla_g p|^2}{p}
\\ \\
& \leq   \int_0^{V_F} \int_0^\infty \big( g_L +g \big) p(v,g,t)   + \frac{1}{\sigma_E} \int_{0}^{+\infty}\int_0^{V_F}   \big( g_{\rm in} - g \big)  \partial_gp \; dgdv.
 \eea $$
 
 We may integrate by parts the last term and find $ \frac{1}{\sigma_E} -  g_{\rm in}  \int_0^{V_F} p(v,0,t)dv$.
Then, using  estimate (\ref{estthmfupper}) of Theorem~\ref{thm_upper}, we deduce that 
$$ \bea
\int_{0}^{V_F} \int_{0}^{+\infty} \ln(p) p(v,g,T) dvdg + \int_{0}^{T} \int_{ g_{F}} ^\infty N(g,t) \ln \f{g V_E }{ gV_E-gV_F- g_L V_F }  + \f 1  {\sg_E}  \int_{0}^{T} \int_{0}^{V_F}\int_{0}^{+\infty}  a(t) \f{|\nabla_g p|^2}{p} 
\\ \\
+ \int_0^T  \int_0^{V_F}  g_{\rm in} (t) p(v,0,t)dv  \leq \int_{0}^{V_F} \int_{0}^{+\infty} \ln(p) p(v,g,0) dvdg +C(T).
\eea$$
A standard argument allows us to recover the absolute value in the logarithm. 
We write
$$ \bea
\int_{0}^{V_F} \int_{0}^{+\infty} |\ln(p)| p(v,g,t) dvdg &=  \int p\ln(p) dvdg -2  \int_{p \leq 1} p \ln(p) dvdg.
\\ \\
& \leq \int p\ln(p) dvdg +  C   \int_{p \leq 1}  p^{\frac{1}{2}} .
\eea $$
It remains to notice that
$$
 \int_{0}^{V_F} \int_{0}^{+\infty} p^{\frac{1}{2}}  dvdg \leq  \int_{0}^{V_F} \int_{0}^{+\infty} (1+g)^2 p(v,g,t)dvdg +  \int_{0}^{V_F} \int_{0}^{+\infty}  (1+g)^{-2} dvdg
$$
and to apply the property of propagation of moments of Theorem~\ref{thm_upper}. This concludes the bound of Theorem~\ref{thm:regularity}.
\hfill $\square$
 
%------------------------------------------------------------------------------
\subsection{Gain of integrability for the total firing}
%----------------------------------------------------------------------------------

We can improve Theorem~\ref{thm_upper} and obtain bounds on $L^q$ norms, $q \geq 2$. More precisely, we prove  that, if the initial data belongs to  $L^q$ for $q \geq 2$, then so does the solution. Moreover the propagation of moments holds in $L^q$. This allows us to control the firing rate  in $L^q, q \geq 2$ assuming the initial data sufficiently decreasing at infinity.
This is stated in the following theorem 
%-------------------------------------------------
\begin{thm}\label{thml2}
Let $\ell \geq 0$, $q \geq 2$,   and assume that the initial data is such that
$$ 
 \int_{0}^{V_F}\int_{0}^{+\infty} (1+g)^{\ell+q-1} p^q(v,g,0)dvdg<+\infty \quad \hbox{ and } \quad   \int_{0}^{V_F}\int_{0}^{+\infty} (1+g)^2 p(v,g,0) dvdg<+\infty.$$
Then, for all $T>0$, the following a priori estimates hold
\begin{equation}\label{est:plq}
 \sup_{t \in [0,T]} \int_{0}^{V_F}\int_{0}^{+\infty} (1+g)^{\ell+q-1} p^q(v,g,t)dvdg <+\infty,
\end{equation}
\begin{equation}\label{est:derplq}
 \int_{0}^{T} \int_{0}^{V_F}\int_{0}^{+\infty} a(t) (1+ g^\ell(-g_Lv+g(V_E-v))_+)^{q-1} (\partial_g p)^2 (v,g,t)p^{q-2} dvdg dt<+\infty.
\end{equation}
  Assume that $\ell > 1+ \frac{q}{q'}$ where $\frac{1}{q}+\frac{1}{q'}=1.$ Then, for all $T>0$,
    \begin{equation}\label{est:Nlq}
 \int_{0}^{T} \mathcal{N}^q(t)dt<+\infty.
\end{equation}
\end{thm}

\noindent{\bf Proof of Theorem~\ref{thml2}.}
We begin with the estimates \eqref{est:plq} and \eqref{est:derplq}. Let  $q \geq 2$ and $\ell \geq 0$. We define the function
$$
K(v,g):= 1 +  g^\ell (-g_Lv+g(V_E-v))_+^{q-1}.
$$
We multiply  equation \eqref{eq:iftdp} by $K(v,g)p^{q-1}$, integrate in $v$ and $g$, use that
$\partial_gK(v,0)=0$, 
to find after integration by parts that 
\begin{eqnarray}\label{eq1}
 \frac{1}{q}\frac{d}{dt} &\int_{0}^{V_F}\int_{0}^{+\infty} K(v,g)p^{q}(v,g,t) dvdg  = - \frac{1}{q} \int_{0}^{+\infty} N(t,g) (p^{q-1}(t,V_F,g)- p^{q-1}(t,0,g))dg  \nonumber 
\\
&+ \int_{0}^{V_F}\int_{0}^{+\infty} R(v,g,t)p^{q}(v,g,t) dvdg  -  (q-1)\frac{ a}{\sigma_E}  \int_{0}^{V_F}\int_{0}^{+\infty} K(v,g)(\partial_g p)^2 p^{q-2} dvdg  
\\ 
&-\frac{a}{\sigma_E} \int_{0}^{V_F}\int_{0}^{+\infty}   \partial_gp \partial_gKp^{q-1}dvdg   -  \frac{(q-1)g_{in}}{q\sigma_E}  \int_{0}^{V_F}K(v,0)p^{q} (t,v,0)dv \nonumber
\end{eqnarray}
where 
$$
R(v,g,t):=  \frac{q-1}{q}(g_L+g)+ \frac{ (-g+g_{in})}{\sigma_E} \partial_gK  + \frac{1}{\sigma_E q}\partial_g[(-g+g_{in})K] .$$
Using the boundary conditions on $p$,  we obtain that 
$$
-\int_{0}^{+\infty} N(t,g) (p^{q-1}(t,V_F,g)- p^{q-1}(t,0,g)) \leq 0.$$
Moreover, as $q \geq 2$, there exists a constant $C$ such that 
$$
R(v,g,t) \leq C(1+g_{in}(t)) (1+g)^{q-1+\ell}.$$
To control the term
$$
 \int_{0}^{V_F}\int_{0}^{+\infty}  a \partial_gp \partial_gKp^{q-1}dvdg,
 $$
we apply the inequality $ab \leq  \varepsilon a^{2} + \frac{1}{\varepsilon} b^2$with  $a=  \partial_gp p^{\frac{q-2}{2}}$, $b= p^{\frac{q}{2}}$ with $\varepsilon >0$ small small enough such that 
$$
\varepsilon \partial_gK <  \frac{q-1}{2} K.
$$
We deduce that, for some constant $C$, we have
$$
\frac{d}{dt} \int_{0}^{V_F}\int_{0}^{+\infty}K(v,g)p^{q}(v,g,t) dvdg \leq C(1+g_{in}(t)) \int_{0}^{V_F}\int_{0}^{+\infty}(1+g)^{q-1+\ell}p^{q}(v,g,t) dvdg$$
$$
- \frac{a(q-1)}{2 \sigma_E}  \int_{0}^{V_F}\int_{0}^{+\infty} K(v,g)(\partial_g p)^2 p^{q-2} dvdg.
$$
 But there exist two positive  constants $C_1, C_2$ such that 
$$ C_1 (1+g)^{q-1+\ell}\leq K(v,g)  \leq C_2 (1+g)^{q-1+\ell}.$$
  Using that $g_{in} \in L^1_{loc}$ (see estimate \eqref{esthk} of Theorem~\ref{thm_upper}), we deduce estimate \eqref{est:plq} and  \eqref{est:derplq} of  Theorem~\ref{thml2}. 
\\  
  
\noindent Let us now prove estimate \eqref{est:Nlq}.  This is a direct consequence of the following Lemma
\begin{lemma}\label{lem:Nq}
  Let $\ell \geq 0$ and $q  \geq 2$. Assume that the initial data is such that 
$$  \int_{0}^{V_F}\int_{0}^{+\infty} (1+g)^{\ell+q} p^q(0,v,g)dvdg<+\infty \hbox{ and }  \int_{0}^{V_F}\int_{0}^{+\infty} (1+g)^2 p(0,v,g)dvdg<+\infty.$$
Then,
$$\int_{0}^{T} \int_{0}^{+\infty} (1+g)^{\ell} N(t,g)^q dgdt<+\infty.$$
\end{lemma}
  Indeed, let us assume for the moment that Lemma~\ref{lem:Nq} holds and let us finish the proof of Theorem~\ref{thml2}.
  We have, using H\"older inequalities, for all $\alpha$ such that $\alpha q'>1$, that
$$ \int_{0}^{T}  \left(\int_{0}^{+\infty} N(t,g) dg \right)^{q} dt \leq  \int_{0}^{T}\left( \int_{0}^{+\infty} (1+g)^{\alpha q}N^q(t,g) dg \right) \left(  \frac{1}{(1+g)^{\alpha q'} }dg\right)^{\frac{q}{q'}} dt$$
and so 
$$ \int_{0}^{T}  \left(\int_{0}^{+\infty} N(t,g) dg \right)^{q} dt \leq C  \int_{0}^{T}\int_{0}^{+\infty} (1+g)^{\alpha q}N^q(t,g) dgdt.$$
Taking $\ell> \frac{q}{q'}$ and using Lemma~\ref{lem:Nq}, we obtain Theorem~\ref{thml2}.  \hfill $\square$

\medskip
  
 \noindent  {\bf Proof of Lemma~\ref{lem:Nq}.}
  We multiply equation  \eqref{eq:iftdp} by $G(v,g)p^{q-1}$
  where
  $$G(v,g):= v(1+g)^{\ell}(-g_Lv+g(V_E-v))_+^{q-1}$$
   to find after integration by parts that there exists a constant $C$ such that
    $$\frac{d}{dt} \int_{0}^{V_F}\int_{0}^{+\infty} G(v,g)p^qdvdg \leq   - \frac{V_F}{q}  \int_{0}^{+\infty}(1+g)^{\ell}N^q(t,g)dg +C(1+g_{in}(t)) \int_{0}^{V_F}\int_{0}^{+\infty}(1+g)^{q+\ell}p^qdvdg $$
    $$ -  \int_{0}^{V_F}\int_{0}^{+\infty} \partial_gG \partial_gp p^{q-1} dvdg. $$
Integrate the above equation in time and using estimates \eqref{est:plq} and  \eqref{est:derplq}, we  obtain Lemma~
\ref{lem:Nq}. 
\hfill $\square$

%%%%%%%%%%%%%%%%%%%%%%%%%%%%%%%%
%---------------------------------------------------------------------------------------------
%\appendix
%---------------------------------------------------------------------------------------------
%%%%%%%%%%%%%%%%%%%%%%%%%%%%%%%%

%--------------------
\bigskip

\noindent {\em Acknowledgment.} The authors would like to thank Albert Cohen for a decisive input in the idea and proof of the Besov estimates of section \ref{sec:stst}.

%%%%%%%%%%%%%%%%%%%%%%%%%%%%%%%%%%%
%
%%%%%% BIBLIO %%%%%%%%%%%%%%%%%%%%%%
%
%%%%%%%%%%%%%%%%%%%%%%%%%%%%%%%%%%%%
%\pagestyle{myheadings}

\end{document}